\documentclass[a4paper, titlepage, reqno, 
11pt
]{proc-l}

\usepackage[T1]{fontenc}
\usepackage[utf8]{inputenc}
\usepackage{lmodern}
\usepackage{xspace}

\usepackage{amssymb}
\usepackage{amsmath,amssymb}
\usepackage{paralist}
\usepackage{graphics}

\usepackage[varqu,varl]{inconsolata}
\usepackage[theoremfont]{newpx}
\usepackage{microtype}

\usepackage{hyperref}
\hypersetup{urlcolor = red, citecolor = blue}
\usepackage[top = 1in,bottom = 1in,left = 1in,right = 1in]{geometry}
\usepackage[sort]{cite}

\usepackage[pagewise]{lineno}
\allowdisplaybreaks[4]
\theoremstyle{plain}
\newtheorem{theorem}{Theorem}[section]
\newtheorem{corollary}[theorem]{Corollary}

\newtheorem{lemma}[theorem]{Lemma}
\newtheorem{proposition}[theorem]{Proposition}

\theoremstyle{definition}

\theoremstyle{remark}
\newtheorem{remark}{Remark}[section]
\numberwithin{equation}{section}

\providecommand{\e}{\ensuremath{\mathrm{e}}}

\newcommand{\uw}{\underline{W}}
\newcommand{\ow}{\overline{W}}
\newcommand{\om}{\overline{M}}
\newcommand{\um}{\underline{M}}
\newcommand{\dd}{\;\mathrm{d}\xspace}

\DeclareSymbolFont{esint}{U}{esint}{m}{n}
\DeclareMathSymbol{\fintsymbol}{\mathop}{esint}{'037}
\def\fint{\fintsymbol}

\newcommand{\loc}{\mathrm{loc}}

\newcommand{\Ma}{M^{(a)}}
\newcommand{\Ra}{R^{(a)}}
\newcommand{\pa}{\psi^{(a)}}
\newcommand{\Ta}{T^{(a)}}

\newcommand{\Mb}{M^{\{b\}}}
\newcommand{\Rb}{R^{\{b\}}}
\newcommand{\pb}{\psi^{\{b\}}}

\begin{document}

\title{A note on the $8\pi$ problem of J\"ager-Luckhaus system}

\author[X.~Mao]{Xuan Mao}
\address[X.~Mao]{School of Mathematics, Hohai University, Nanjing 211100, Jiangsu, China}
\address[X.~Mao]{Laboratory of Mathematical Modeling and Intelligent Computing for Water Systems, Hohai University, Nanjing 211100, China}
\email[X.~Mao]{maox@hhu.edu.cn}
\urladdr[X.~Mao]{https://orcid.org/0009-0001-4960-8901}
%\thanks{$^\ast$ Corresponding author.}
\thanks{Xuan~Mao is funded by 
``the Fundamental Research Funds for the Central Universities'' (No.~B250201215)
and Basic Research Program of Jiangsu~(No.~BK20251482).}

\author[M.~Liu]{Meng Liu}
\address[M.~Liu]{Department of Applied Mathematics, Anhui University of Technology, Ma'anshan
243002, P. R. China}
\email[M.~Liu]{LMeng\_math@yeah.net}
\thanks{Meng Liu has been supported by the Youth Research Fund of Anhui University of Technology (No. QZ202421, No. QD202377 ).}

\author[Y.~Li]{Yuxiang Li}
\address[Y.~Li]{School of Mathematics, Southeast University, Nanjing 211189, P. R. China}
\email[Y.~Li]{lieyx@seu.edu.cn}
\thanks{Yuxiang Li is supported in part by National Natural Science Foundation of China (Nos. 12271092 and 11671079) 
and the Jiangsu Provincial Scientific Research Center of Applied Mathematics under Grant No. BK20233002.}


% 35B40 Asymptotic behavior of solutions to PDEs 
% 35B33 Critical exponents in context of PDEs
% 35B35 Stability in context of PDEs
% 35B44 Blow-up in context of PDEs
% 35K65 Degenerate parabolic equations
% 92C17 Cell movement (chemotaxis, etc.)

% \date{}%
% \dedicatory{}%
% \commby{}%
% ----------------------------------------------------------------

\begin{abstract}
We show that for any nonnegative, radially symmetric and continuous initial datum with critical mass $8\pi$, 
J\"ager-Luckhaus system in the unit disk, 
known as a parabolic-elliptic Keller-Segel model, 
admits a globally bounded classical solution. 
Moreover, it is asserted that the spatial constant equilibrium $8$ is globally and exponentially asymptotically stable.
\end{abstract}

\maketitle

\section{Introduction}
\label{introduce section}

This paper is concerned with J\"ager-Luckhaus system 
\begin{align}
\begin{cases}
	\label{sys: JL}
	u_t = \Delta u - \nabla \cdot(u\nabla v),&  x\in\Omega, t>0, \\
	0 =  \Delta v - \mu + u,&  x\in\Omega,	t>0, \\
	\partial_\nu u = \partial_\nu v = 0 , & x\in\partial\Omega, t >0, \\
	u(\cdot, 0) = u_0, & x\in\Omega,
\end{cases}
\end{align}
where $\mu := \fint_\Omega u_0\dd{x}$ denotes the average initial mass,
$\nu$ denotes the outward unit normal vector field to the boundary $\partial\Omega$,
and $\Omega\subset\mathbb{R}^2$ is the unit disk $B_1:=\{x\in\mathbb{R}^2: |x|<1\}$. 
The system \eqref{sys: JL} is a parabolic-elliptic simplified version of Keller-Segel model~\cite{Keller1970}, proposed by J\"ager and Luckhaus~\cite{Jaeger1992}. 
It describes oriented movement of cells in response to chemical stimulus.
Here, $u = u(x, t)$ denotes the density of cells and $v = v(x, t)$ represents the concentration of chemical signals. 
The goal of this paper is to investigate the global classical solvability of \eqref{sys: JL} with critical initial mass $8\pi$.

%%%%%%%%%%%%%%%%%%%%%%%%%%%%%%%%%%%%%%%%%%%%%%%%%%%%%%%%%%%%%%%%%%%%%%%%%%%%%%%%%%%%

J\"ager and Luckhaus \cite{Jaeger1992} found that radially symmetric initial function with large mass can be constructed such that the corresponding solution explodes at $x=0$ in finite time, 
whereas for a large class of continuous initial data $u_0$ with mass less than a number, 
solutions of \eqref{sys: JL} exist globally and remain bounded.
Nagai \cite{Nagai1995} considered radial solutions of \eqref{sys: JL} and deduced the two-dimensional critical mass $8\pi$ for finite-time blow-up in the sense that 
in the $n$-dimensional open ball $\Omega = \{x\in\mathbb{R}^n: |x|<R\}$,
\begin{itemize}
\item either $n=1$ or $n=2$ and $\int_\Omega u_0\dd{x}<8\pi$ implies all solutions are globally bounded;
\item when $n\geq3$ and $m>0$ or $n=2$ and $m>8\pi$, there exists $\varepsilon(m)>0$ such that the solution of \eqref{sys: JL} with initial data satisfying $\int_\Omega u_0\dd{x} = m$ and $\int_\Omega u_0|x|^2\dd{x} <\varepsilon(m)$ blows up in finite time.
\end{itemize}
Ohtsuka, Senba and Suzuki \cite{Ohtsuka2007} showed that 
if a radial solution with supercritical mass blows up, 
then blow-up happens in finite time.
Winkler considered instability of spatial homogeneity 
and deduced that for $n\ge2$, there exists $m_\star(n,R)>0$ such that 
for any radial initial datum with mass exceeding $m_\star$,
which is more concentrated than its average mass in an appropriate sense, 
the solution blows up in finite time~\cite{Winkler2019}.
Souplet and Winkler studied the spatial profile of radial solutions exploding in finite time, under mild assumptions on initial data in higher dimensions~\cite{Souplet2019}.
For more related results, 
we refer to \cite{Nagai2001, Nagai1997, Nagai2000, Senba2000,Horstmann2001, Winkler2013, Kohatsu2025} and references therein. 

%%%%%%%%%%%%%%%%%%%%%%%%%%%%%%%%%%%%%%%%%%%%%%%%%%%%%%%%%%%%%%%%%%%%%%%%%%%

To the best of our knowledge,
only two special parabolic-elliptic variants of Keller-Segel system have available information on dynamic behaviors 
when the initial mass is critical.
When $\Omega = \mathbb{R}^2$, 
the system \eqref{sys: JL} is reduced to the Patlak-Keller-Segel model
\begin{equation}\label{sys: pks}
  \begin{cases}
    u_t = \Delta u - \nabla\cdot(u\nabla v),  & x\in\mathbb{R}^2, t>0,\\
    v := -\frac{1}{2\pi}\ln |x| * u, & x\in\mathbb{R}^2, t>0.
  \end{cases}
\end{equation}
$8\pi$ is the critical mass for global well-posedness of mild solutions to the system \eqref{sys: pks} \cite{Wei2018,Blanchet2006}. 
The system \eqref{sys: pks} has a family of stationary solutions \cite{Chen1991}
\begin{equation}
  \label{sym: special solution of pks}
  U_{\lambda, x_0}(x) = \frac1{\lambda^2}U_0\bigg(\frac{x-x_0}{\lambda}\bigg), \quad U_0(x) := \frac{8}{(1+|x|^2)^2}, 
  \quad \lambda > 0, x_0\in\mathbb{R}^2.
\end{equation}
There exist initial data with critical mass such that the solution approaches to certain stationary solutions \cite{Blanchet2012,Carlen2013,Lopez-Gomez2013}, 
as time goes to infinity,  
that chemotaxis collapse, 
blow-up in the form of Dirac-delta distribution with $8\pi$ weight, 
happens in infinite time at the center of mass \cite{Blanchet2008}, 
and such that the solution oscillates among stationary solutions \cite{Lopez-Gomez2014}. Refined descriptions of infinite-time blow-up were given in \cite{Ghoul2018,Davila2024}.
When $\Omega$ is the unit disk, the system \eqref{sys: pks} is known as Smoluchowski-Poisson system
\begin{equation}
  \label{sys: ks-SP}
  \begin{cases}
  u_t = \Delta u - \nabla\cdot(u\nabla v), & x\in\Omega, t>0,\\
  0 = \Delta v + u, & x\in\Omega, t>0,
  \end{cases}
\end{equation} 
subjected to mixed boundary conditions
\begin{equation*}
  \partial_\nu u - u\partial_\nu v = v = 0
  \quad\text{for } x\in\partial\Omega \text{ and } t>0.
\end{equation*} 
Biler, Karch, Lauren\c{c}ot and Nadzieja showed that 
the radial solutions of the system~\eqref{sys: ks-SP} exists globally 
and a chemotactic collapse takes place in infinite time~\cite{Biler2006a}. 
When comparing with companion work~\cite{Biler2006}, they predicted that
the solutions of the system~\eqref{sys: ks-SP} should grow up faster than that of~\eqref{sys: pks}. 
A rate of grow-up of~\eqref{sys: ks-SP} was given by Kavallaris and Souplet~\cite{Kavallaris2009},
which is indeed faster than that of~\eqref{sys: pks} given by Ghoul and Masmoudi~\cite{Ghoul2018}.
Without symmetry assumptions, Suzuki~\cite{Suzuki2013} ruled out possibility of boundary blow-up and showed that $8\pi$ is the critical mass of \eqref{sys: ks-SP} for finite-time blow-up, 
and that if initial mass is $8\pi$, 
then the solution exhibits a collapse in infinite time,
 of which movement is subjected to a gradient flow associated with the Robin function. 
%%%%%%%%%%%%%%%%%%%%%%%%%%%%%%%%%%%%%%%%%%%%%%%%%%%%%%%%%%%%%%%%%

Previous works suggest that when initial mass is critical, 
the system \eqref{sys: JL} should admit a collapse in infinite time.
However, in this paper, 
we found that under the framework of radially symmetric solutions, 
the system \eqref{sys: JL} with critical mass has no collapse at all.

\begin{theorem}
  \label{thm: global wellposedness}
  Let $\Omega=B_1\subset\mathbb{R}^2$. 
  Suppose that the initial function $u_0\in C^0(\overline{\Omega})$ is nonnegative and radially symmetric such that 
  \begin{equation}
  \int_\Omega u_0\dd x \in(0, 8\pi].
  \end{equation}
  Then the classical solution of the system~\eqref{sys: JL}, given by Proposition~\ref{prop: local well-posedness} below, exists globally and remains bounded.
  Moreover, it converges to the unique constant steady state exponentially, i.e., 
  \begin{equation}
    \label{eq: exponential stability}
    \sup_{t>0}\bigg\{\e^{\ell t}\bigg\Vert u-\fint_\Omega u_0\dd x \bigg\Vert_{L^\infty(\Omega)} \bigg\}
    < \infty,
  \end{equation}
  for certain constant $\ell > 0$ independent of initial data $u_0\in C^0(\overline{\Omega})$.
\end{theorem}

\begin{remark}
Inspired by precedents~\cite{Jaeger1992,Biler2006a}, 
we study the Dirichlet initial boundary value problem 
\begin{equation*}
    \begin{cases}
    M_t = 4\xi M_{\xi\xi} 
    + \frac{MM_\xi}{\pi} 
    - \frac{m\xi M_\xi}{\pi}, & (\xi,t) \in (0,1)\times(0,T_{\max}),\\
    M(0,t) = 0, M(1,t) = m, & t \in (0,T_{\max}),\\
    M(\cdot, 0) = M_0(\cdot), & \xi \in (0,1),
    \end{cases}
  \end{equation*}
satisfied by the mass distribution function
\begin{equation*}
  M(\xi,t) := \int_{B_{\sqrt{\xi}}}u(x,t)\dd x
  = 2\pi\int_0^{\sqrt\xi} u(r,t)r\dd r
  \quad\text{for } (\xi,t) \in [0,1]\times[0,T_{\max}),
\end{equation*}
where $u$ is the component of a classical solution to the system~\eqref{sys: JL}.
We construct a family of stationary supersolutions to the scalar parabolic problem within the framework of (sub)critical mass.
This is sufficient to show global boundedness via an $\varepsilon$-regularity.
Combining a family of stationary subsolutions,
we are able to present an alternative proof to uniqueness of radially symmetric stationary solutions of the system \eqref{sys: JL} with (sub)critical mass~\cite[Theorem~2.2]{Horstmann2011}.
It follows from LaSalle invariant principle that the semi-trivial steady state $m/\pi$ for $m\in(0,8\pi]$ is globally asymptotically stable.
We further obtain that the solutions emanating from both supersolutions and subsolutions approach to $m\xi$ exponentially. 
To establish exponential convergence for generic radially symmetric initial data, 
it suffices to employ appropriate comparison and interpolation arguments.
\end{remark}

\begin{remark}
  Our comparison-based methods, while relying heavily on radial symmetry, 
  achieve delicate results beyond global solvability,
  especially when contrasted with related theories \cite{Suzuki2025} valid to general situations. 
  % For the global-in-time result of ~\eqref{sys: ks-SP}, one assumes the blowup of the solution
  % in finite time to get the formation of collapse with the quantized mass, $8\pi$
  % times an integer, with the positive residual term, which ensures $\|u_0\|_1>8\pi$, a contradiction. This argument is valid to more general situations, without
  %  radial symmetry or in the other form of Possion equation, see \cite{Suzuki2025} for more details.
\end{remark}

This paper is organized as follows. 
In Section \ref{sec: preliminary}, 
we introduce some preliminaries. 
Section~\ref{sec: exclusion of blowup} is devoted to global boundedness of solutions to the system \eqref{sys: JL}. 
We present an alternative proof of uniqueness of radially symmetric stationary solutions to the system \eqref{sys: JL} with (sub)critical mass
in Section~\ref{sec: uniquenss and stability of ss}. 
We finish the proof of 
Theorem~\ref{thm: global wellposedness} by offering the exponential convergence.

\section{Preliminaries}
\label{sec: preliminary}
We shall first state results on local existence and uniqueness of classical solutions to the system \eqref{sys: JL}, and refer to \cite{Cieslak2008,Winkler2010b} for complete proofs.

\begin{proposition}
  \label{prop: local well-posedness}
Let $\Omega=B_1\subset\mathbb{R}^2$. 
Suppose that $u_0\in C^0(\overline{\Omega})$  is nonnegative and radially symmetric. 
Then there exist
\(T_{\max}\in(0,\infty]\) and a pair $(u, v)$ of radially symmetric functions solving \eqref{sys: JL} in the classical sense, 
such that 
$(u, v)$ can be uniquely identified by the inclusions		
	\begin{align*}
	& u \in C^0(\overline{\Omega} \times[0, T_{\max})) \cap C^{2,1}(\overline{\Omega} \times(0, T_{\max})), \\
	& v \in \cap_{q>2}L^\infty_{\loc}([0,T_{\max}),W^{1,q}(\Omega))\cap C^{2,0}(\overline{\Omega} \times(0, T_{\max})), 
	\end{align*}
and the identity
\begin{equation*}
\int_{\Omega} v(\cdot, t)=0 
\quad\text {for all } t \in[0,T_{\max}).
\end{equation*}
Furthermore, \(u>0\) in \(\overline{\Omega}\times(0,T_{\max})\),
\begin{equation*}
	\int_{\Omega} u(x, t) \dd {x}=\int_{\Omega} u_{0} \dd x \quad\text {for all } t \in (0,T_{\max}),
\end{equation*}
and 
\begin{equation*}
\text{if } T_{\max}<\infty, \quad\text{then } \limsup_{t\nearrow T_{\max}}\|u\|_{L^\infty(\Omega)} = \infty.
\end{equation*}
\end{proposition}

Following precedents~\cite{Jaeger1992,Biler2006}, 
we shall analyze a scalar parabolic problem satisfied by the mass distribution function~\eqref{sym: mass distribution function}. 

\begin{lemma}
  \label{le: mass distribution function}
  Let $\Omega = B_1\subset\mathbb{R}^2$ and $(u,v)$ be the radially symmetric classical solution given by Proposition~\ref{prop: local well-posedness}.
  Then the mass distribution function
  \begin{equation}
    \label{sym: mass distribution function}
    M(\xi,t) := \int_{B_{\sqrt{\xi}}}u(x,t)\dd x
    = 2\pi\int_0^{\sqrt\xi} u(r,t)r\dd r
    \quad\text{for } (\xi,t)\in[0,1]\times[0,T_{\max}),
  \end{equation}
  solves the Dirichlet initial boundary value problem
\begin{equation}
\label{sys: scalar system derived from JL system}
  \begin{cases}
  M_t = 4\xi M_{\xi\xi} 
  + \frac{MM_\xi}{\pi} 
  - \frac{m\xi M_\xi}{\pi}, &(\xi,t) \in (0,1)\times(0,T_{\max}),\\
  M(0,t) = 0, M(1,t) = m, & t\in(0,T_{\max}),\\
  M(\cdot, 0) = M_0(\cdot), & \xi\in(0,1),
  \end{cases}
\end{equation}
where 
\begin{equation*} 
m := \int_\Omega u_0\dd x
\quad\text{and}\quad
M_0(\xi) := \int_{B_{\sqrt\xi}}u_0(x)\dd x
= 2\pi\int_0^{\sqrt{\xi}}u_0(r)r\dd r.
\end{equation*}
\end{lemma}

The degenerate parabolic problem \eqref{sys: scalar system derived from JL system} admits a comparison argument 
(see e.g., \cite[Lemma~5.1]{Bellomo2017}).

\begin{lemma}\label{le: comparision principle}
Let $T\in(0,\infty]$ and suppose that $\um$ and $\om$ are two nonnegative functions 
which belong to 
\(C^1([0, 1]\times [0, T))\cap C^2((0, 1)\times(0, T))\) 
and satisfy
\begin{equation*}
\um_t - 4\xi \um_{\xi\xi} - \frac{\um\um_\xi}{\pi} + \frac{m\xi\um_\xi}{\pi} 
\leq 
\om_t - 4\xi \om_{\xi\xi} - \frac{\om\om_\xi}{\pi} + \frac{m\xi\om_\xi}{\pi}  
\quad\text{in }(0,1) \times (0, T),
\end{equation*}
and if moreover
\begin{equation*}
	\um(0, t) \leq \om(0, t) 
\quad{and}\quad
\um(1, t) \leq \om(1, t) 
\quad\text{in } (0, T),
\end{equation*}
as well as
\begin{equation*}
\um(\xi, 0)\leq \om(\xi, 0) 
\quad\text{in }  (0, 1),
\end{equation*}
then
\begin{equation*}
\um (\xi, t) \leq \om (\xi, t) 
\quad\text{in } [0, 1]\times [0, T).
\end{equation*}
\end{lemma}

Analogous to the classical Keller-Segel system~\cite{Nagai1997},
the J\"ager-Luckhaus system~\eqref{sys: JL} has a Lyapunov functional.  
We shall give a direct deduction for convenience.
\begin{lemma}
\label{le: free energy functional}
Let \(\Omega = B_1\subset\mathbb{R}^2\).
Then for the solution \((u,v)\) given by Proposition~\ref{prop: local well-posedness},
\begin{equation}
\label{eq: differential equation of free energy}
  \frac{\dd}{\dd t} \mathcal{F}(u)
  =-\int_\Omega u|\nabla(\ln u - v)|^2
\end{equation}
holds for all  $t\in(0,T_{\max})$,
where 
\begin{equation}
  \label{sym: free energy}
  \mathcal{F}(u) := \int_\Omega u\ln u 
  - \frac{1}{2}\int_\Omega uv
  \quad \text{for } u\in C^0(\overline{\Omega}).
\end{equation}
Here and below, $v:=(-\Delta )^{-1} (u-\mu)\in\{f\in W^{2,2}(\Omega)\mid \int_\Omega f = 0 \text{ and } \partial_\nu f = 0 \text{ on }\partial\Omega \}$ 
denotes the unique solution of Poisson equation $-\Delta v = u-\mu$ with homogeneous Neumann boundary conditions.  
\end{lemma}

\begin{proof}
  Multiplying the first equation in \eqref{sys: JL} by \(\ln u - v\) and integrating by parts,
  we get
  \begin{equation}
  \label{eq: differentiate free energy}
    \begin{aligned}
    \frac{\dd}{\dd t}\int_\Omega (u\ln u - uv) 
    &= \int_\Omega u_t(\ln u - v) + \int_\Omega u_t - \int_\Omega uv_t\\
    &= -\int_\Omega u|\nabla (\ln u - v)|^2 - \int_\Omega uv_t.
    \end{aligned}
  \end{equation}
  Using the second equation in \eqref{sys: JL}, we have
  \begin{equation*}
    \int_\Omega uv_t 
    = \int_\Omega (-\Delta v + \mu) v_t
    = \frac12\frac{\dd}{\dd t}\int_\Omega|\nabla v|^2.
  \end{equation*}
  Noting
  \begin{equation*}
    0
    =\int_\Omega(\Delta v - \mu + u)v
    = \int_\Omega uv - \int_\Omega |\nabla v|^2, 
  \end{equation*}
  we deduce \eqref{eq: differential equation of free energy} by substituting two preceding identities into \eqref{eq: differentiate free energy}. 
\end{proof}

\section{Global boundedness}
\label{sec: exclusion of blowup}

This section is devoted to the global boundedness of solutions to the system \eqref{sys: JL} with (sub)critical mass \(m\in(0,8\pi]\).

\begin{proposition}
  \label{prop: global boundedness}
  Let $\Omega = B_1\subset\mathbb{R}^2$. 
  Suppose that the initial function $u_0\in C^0(\overline{\Omega})$ is radially symmetric and nonnegative such that 
\begin{equation*}
  \int_\Omega u_0\dd x = m \in(0, 8\pi].
\end{equation*} 
Then the unique solution $(u, v)$ given by Proposition~\ref{prop: local well-posedness} exists globally and remains bounded.
\end{proposition}

To this end, we construct a family of supersolutions  
\begin{equation}
  \label{sym: supersolution of steady states problem} 
  \ow_a(\xi) 
  := \frac{m(a+1)\xi}{a+\xi}
  \quad\text{for } a > 0 \text{ and } \xi\in(0,1), 
\end{equation}
to the stationary problem
\begin{equation}
\label{sys: stationary problem of scalar parabolic problem}
  \begin{cases}
    \mathcal{Q}W = 0, & \xi\in (0,1), \\
    W(0) = 0,\quad W(1) = m,
  \end{cases}
\end{equation}
in the sense of \eqref{eq: ow is a strict stationary supersolution} below,
where
\begin{equation}\label{sym: operator Q}
  \mathcal{Q}M := 
  - 4\xi M_{\xi\xi} 
  - \frac{MM_\xi}{\pi} 
  + \frac{m\xi M_\xi}\pi.
\end{equation}
Using some supersolution \eqref{sym: supersolution of steady states problem} to bound from above initial data \(M_0\) (up to a later time),
we then can use comparison method (see Lemma~\ref{le: comparision principle}) and \(\varepsilon\)-regularity (see Lemma~\ref{le: boundedness principle}) to establish global well-posedness and boundedness.
Noting the mass distribution functions of exact solutions~\eqref{sym: special solution of pks}  with $x_0=0$, are given by 
\begin{equation}
  \label{eq: accumulated density of special solutions}
  \frac{8\pi \xi}{\lambda^2 + \xi}
  \quad \text{for } \xi > 0 \text{ and } \lambda > 0,
\end{equation}
we are inspired to construct~\eqref{sym: supersolution of steady states problem}, 
as a perturbation of special functions~\eqref{eq: accumulated density of special solutions}. 

Our key ingredient is the following observation.
\begin{lemma}
\label{le: a family of supersolutions of steady states}
Let \(m\in(0,8\pi]\). 
For each \(a>0\), 
\(\ow_a\) is increasing and concave in \((0,1)\),
which satisfies
\begin{equation}
  \label{eq: ow asymptotic}
  \lim_{a\searrow0}\ow_a(\xi) = m,
  \quad
  \lim_{a\to\infty}\ow_a(\xi) = m\xi
  \quad\text{for all } \xi\in(0,1)
\end{equation}
and
\begin{equation}
\label{eq: ow is a strict stationary supersolution}
  \mathcal{Q}\ow_a > 0\quad\text{for all } \xi \in (0,1).
\end{equation}
\end{lemma}

\begin{proof}
  \eqref{eq: ow asymptotic} is obvious. 
  It suffices to check \eqref{eq: ow is a strict stationary supersolution}. 
  Noting 
  \begin{equation*}
    \ow_{a\xi} = 
    \frac{m(a+1)}{a+\xi} - \frac{m(a+1)\xi}{(a+\xi)^2}
    = \frac{m(a+1)a}{(a+\xi)^2} 
  \end{equation*}
  and
  \begin{equation*}
    \ow_{a\xi\xi}  
    = -\frac{2m(a+1)a}{(a+\xi)^3}, 
  \end{equation*} 
  we calculate and estimate by \(m\in(0,8\pi]\)
  \begin{align*}
   \mathcal{Q}\ow_a &= 
    \frac{8m\xi(a+1)a}{(a+\xi)^3}
    - \frac{m^2\xi(a+1)^2a}{\pi(a+\xi)^3}
    + \frac{m^2\xi(a+1)a}{\pi(a+\xi)^2}\\
    &= \frac{m\xi(a+1)a}{(a+\xi)^3}
    \cdot\left[8 - \frac{m(a+1)}{\pi} + \frac{m(a+\xi)}{\pi}\right]\\
    &= \frac{m\xi(a+1)a}{(a+\xi)^3}
    \cdot \frac{8\pi- m + m\xi}{\pi}  > 0
    \quad\text{for all } \xi\in(0,1)
  \end{align*} 
  and thereby complete the proof.
\end{proof}

The following auxiliary arguments are used to bound the mass distribution function $M(\xi, T)$ from Lemma~\ref{le: mass distribution function} for some intermediate time $T\in[0,T_{\max})$.

\begin{lemma}
  \label{le: function with bound derivative can be bounded}
  Suppose that $f\in C^0([0,1])\cap C^1((0,1))$ complies with
  \begin{equation*}
    \frac{1}{C} < f' < C\quad\text{for some }C > f(1)-f(0) > 0.  
  \end{equation*} 
  Let the set 
  \begin{equation*}
    \{f_a: [0,1]\mapsto[f(0),f(1)]\mid f_a(0)=f(0), f_a(1)=f(1)\}_{a\in \mathcal{I}}
  \end{equation*} 
  be a family of concave functions satisfying 
  \begin{equation}
  \label{eq: property of a family concove functions}
    \sup_{a\in\mathcal{I}} f_a(\xi) = f(1)\quad\text{for all } \xi\in(0,1).
  \end{equation}
  Then there exists $\alpha\in\mathcal{I}$ such that
  \begin{equation*}
    f(\xi) < f_\alpha(\xi) \quad\text{for all } \xi \in (0,1).  
  \end{equation*}
  Similarly, if the set  
  \begin{equation*} 
    \{f_b:[0,1]\mapsto[f(0),f(1)]\mid f_b(0)=f(0), f_b(1)=f(1)\}_{b\in\mathcal{J}}
  \end{equation*}
  consists of a family of convex functions satisfying
  \begin{equation*}
    \inf_{b\in\mathcal{J}} f_b(\xi) = f(0) \quad\text{for all } \xi\in(0,1),
  \end{equation*}
  then there exists $\beta\in\mathcal{J}$ such that
  $f(\xi) > f_\beta(\xi)$ for all $\xi\in(0,1)$.
\end{lemma}
\begin{proof}
  Using Newton-Leibniz formula, we estimate
  \begin{equation*} 
  f(\xi) = f(0) + \int_0^\xi f'\dd \eta  
  \leq f(0) + C\xi
  \quad\text{for all } \xi\in(0,1),
  \end{equation*}
  and
  \begin{equation*}
    f(\xi) = f(1) - \int_\xi^1 f'\dd \eta 
    \leq f(1) - \frac{1-\xi}{C}
    \quad\text{for all } \xi \in (0,1).
  \end{equation*}
  So we obtain 
  \begin{equation*}
    f\leq \overline{f} := 
    \min\bigg\{f(0)+C\xi, f(1) - \frac{1-\xi}{C}\bigg\}
    \quad\text{for all } \xi\in(0,1).
  \end{equation*}
  We note that
  \begin{equation*} 
  \xi_0 := \frac{f(1)-f(0)-\frac1C}{C-\frac1C}\in(0,1),
  \end{equation*}
  solves the equation
  \begin{equation*} 
  f(0) + C\xi = f(1) - \frac{1-\xi}{C}, \quad\xi\in\mathbb{R},
  \end{equation*}
  and thus $\overline{f}(\xi_0)\in(f(0),f(1))$.
  Thanks to \eqref{eq: property of a family concove functions},
  there exists $\alpha\in\mathcal{I}$ such that 
  $f_\alpha(\xi_0) > \overline{f}(\xi_0)$.
  In light of the concavity of \(f_\alpha\),
  \begin{equation}
    \label{eq: f alpha >} 
  \frac{f_\alpha(\xi)-f(0)}{\xi} 
  \geq \frac{f_\alpha(\xi_0)-f(0)}{\xi_0}
  > \frac{\bar{f}(\xi_0)-f(0)}{\xi_0} = C
  \quad\text{for all } \xi\in (0,\xi_0)
  \end{equation}
  and
  \begin{equation} 
    \label{eq: f alpha <}
  \frac{f_\alpha(1)-f_\alpha(\xi)}{1-\xi}
  \leq \frac{f_\alpha(1)-f_\alpha(\xi_0)}{1-\xi_0}
  < \frac{f(1)-\bar{f}(\xi_0)}{1-\xi_0} 
  = \frac{1}{C}
  \quad\text{for all }\xi\in(\xi_0,1).
  \end{equation}
Rearranging \eqref{eq: f alpha >} and \eqref{eq: f alpha <} gives us  
  \begin{equation*} 
    f_\alpha > \bar{f} \geq f\quad\text{for all } \xi \in (0,1).
  \end{equation*}
  One may conclude the second assertion analogously, 
  so we suppress the proof here.
\end{proof}

\begin{lemma}
  \label{le: a upper bound of initial data}
  Let $\Omega = B_1\subset\mathbb{R}^2$. 
  Suppose that the initial function $u_0\in C^0(\overline{\Omega})$ is radially symmetric and nonnegative such that 
\begin{equation*}
  \int_\Omega u_0\dd x = m \in(0, 8\pi].
\end{equation*} 
Then for each $T\in(0,T_{\max})$,
there exists $\alpha>0$ such that 
\begin{equation}
  \label{eq: a upper bound of initial data}
  M(\xi, T) \leq \ow_\alpha(\xi)\quad\text{for all } \xi\in(0,1),
\end{equation} 
where $M$ and $\ow_\alpha$ are defined in \eqref{sym: mass distribution function} and \eqref{sym: supersolution of steady states problem}, respectively.
\end{lemma}
\begin{proof}
  By the strong maximal principle,
  $u > 0$ holds over $\overline{\Omega}\times(0,T_{\max})$.
  Hence, for any $T\in(0,T_{\max})$, 
  there exists $C>0$ such that 
  \begin{equation*}
   \frac1C < u(x,T) < C,\quad \text{for all } x\in\overline{\Omega}. 
  \end{equation*}
  Since $M(\xi, t)$ given by \eqref{sym: mass distribution function} satisfies
  $M(0,t)=0$, 
  $M(1, t) = m$ 
  and $M_\xi(\xi,t) = \pi u(\sqrt{\xi}, t)$,
  and Lemma~\ref{le: a family of supersolutions of steady states} tells that $\ow_a$ is concave for each $a>0$ and satisfies
  \eqref{eq: ow asymptotic},
  we infer from Lemma~\ref{le: function with bound derivative can be bounded} that there exists $\alpha>0$ such that 
  \eqref{eq: a upper bound of initial data} is valid.
\end{proof}

The following $\varepsilon$-regularity was shown in \cite[Lemma~4.1]{Winkler2019} via a Bernstein-type argument. 
Here, we would like to give a shorter proof inspired by Nagai~\cite{Nagai1995}.
\begin{lemma}
\label{le: boundedness principle}
Let $\Omega=B_1\subset\mathbb{R}^2$, and suppose that $u_0\in C^0(\overline{\Omega})$ is nonnegative and radially symmetric. 
Assume that the function $M$ defined in \eqref{sym: mass distribution function} has the property that
\begin{equation*} 
L := \sup_{(\xi, t) \in(0, 1) \times(0, T_{\max })} 
\frac{M(\xi, t)}{\xi}<\infty.
\end{equation*}
Then there exists a positive constant $C = C(L,m,\|u_0\|_{L^\infty(\Omega)})>0$ such that
\begin{equation*} 
M_\xi(\xi, t) \leq C 
\quad \text {for all } \xi \in(0,1)
\text { and } t \in(0, T_{\max}). 
\end{equation*}
\end{lemma}
\begin{proof}  
  Integrating the second equation $-(v_rr)_r = ur - mr/\pi$ in \eqref{sys: JL} from $r=0$ to $\sqrt\xi$,
  we obtain
  \begin{equation*} 
   - 2\pi v_r(\sqrt\xi, t)\sqrt\xi 
   = 2\pi\int_0^{\sqrt\xi} ur\dd r - 2\pi\int_0^{\sqrt\xi}\frac{mr}{\pi}\dd r
   = M(\xi, t) - m\xi
  \end{equation*}
  for all $(\xi,t)\in(0,1)\times(0,T_{\max})$.
  This implies
  \begin{equation*}
  2\pi|\nabla v(x,t)| = 2\pi|v_r(|x|,t)|
  \leq \frac{M(|x|^2,t)}{|x|^2}\cdot|x| + m|x|
  \leq %\sup_{(\xi, t) \in(0, 1) \times(0, T_{\max })} \frac{M(\xi, t)}{\xi} + m 
  L + m
  \end{equation*}
  for all $(x,t)\in\Omega\times(0,T_{\max})$. 
  By the well-established Neumann heat semigroup estimates~\cite[Lemma~1.3]{Winkler2010a},
  we can find a positive constant $C_{\mathrm{N}} > 0$, independent of $T\in(0,T_{\max})$, such that  
  \begin{align*}
    \|u\|_{L^\infty(\Omega)} 
    &\leq \|\e^{t\Delta} u_0\|_{L^\infty(\Omega)} 
    + \int_0^t\|\e^{(t-s)\Delta}\nabla\cdot(u\nabla v)\|_{L^\infty(\Omega)}\dd s\\
    &\leq \|u_0\|_{L^\infty(\Omega)} 
    + C_{\mathrm{N}}(L+m)\int_0^t(1+s^{-5/6})\e^{-\lambda s}\dd s \cdot\sup_{s\in(0,T)}\|u\|_{L^3(\Omega)}\\
    &\leq \|u_0\|_{L^\infty(\Omega)} 
    + C_{\mathrm{N}}(L+m)m^{1/3}\int_0^\infty(1+s^{-5/6})\e^{-\lambda s}\dd s \cdot\sup_{s\in(0,T)}\|u\|_{L^\infty(\Omega)}^{2/3} 
  \end{align*}
  for all $t\in(0,T)$,
  where $\lambda>0$ is the first nonzero eigenvalue of minus Neumann Laplacian.
  Thus we have 
  \begin{equation*}
    \sup_{s\in(0,T)}\|u\|_{L^\infty(\Omega)} 
    \leq 2\|u_0\|_{L^\infty(\Omega)}+\left(C_{\mathrm{N}}(L+m)m^{1/3}\int_0^\infty(1+s^{-5/6})\e^{-\lambda s}\dd s\right)^3:= \frac{C}{\pi}
  \end{equation*}
  for all $T\in(0,T_{\max})$.
  The desired result follows.
\end{proof}

Now, we are in a position to prove Proposition~\ref{prop: global boundedness}.

\begin{proof}[Proof of \upshape{Proposition~\ref{prop: global boundedness}}]
  By Lemma~\ref{le: a upper bound of initial data},
  there exists \(a>0\) such that
  \begin{equation*}
  M(\xi, T) \leq \ow_a(\xi)
  \quad\text{for all }\xi\in(0,1),
  \end{equation*}
  with 
  \begin{equation*}
  T = \min\left\{1, \frac{T_{\max}}{2}\right\}.
  \end{equation*}
  Recalling Lemma~\ref{le: a family of supersolutions of steady states}, 
  we have by the comparison method stated in Lemma~\ref{le: comparision principle}, 
  \begin{equation*}
  M(\xi, t) \leq \ow_a(\xi) 
  \quad\text{for all } 
  (\xi,t) \in [0,1]\times[T,T_{\max}),
  \end{equation*}
  which implies
  \begin{equation*}
    \begin{aligned}
    \frac{M(\xi,t)}{\xi} 
    \leq \frac{\ow_a(\xi)}{\xi}
    = \frac{m(a+1)}{a+\xi}
    \leq \frac{m(a+1)}{a}  
    \quad\text{for all } 
  (\xi,t) \in [0,1]\times[T,T_{\max}).
    \end{aligned}
  \end{equation*}
  Since $u$ is obviously bounded in $[0,1]\times[0,T]$, 
  Lemma~\ref{le: boundedness principle} entails that there exists \(C>0\) such that
  \begin{equation*} 
    u(x,t) = \frac{M_\xi(|x|^2, t)}{\pi} < C
    \quad\text{for all }(x, t)\in\Omega\times(0,T_{\max}), 
  \end{equation*} 
  which yields $T_{\max} = \infty$.
\end{proof}

\section{Uniqueness and stability of radially symmetric steady states}
\label{sec: uniquenss and stability of ss}

The goal of this section is to present an alternative proof of uniqueness of radial solutions to the stationary problem 
\begin{equation}
  \label{sys: stationary problem of uv}
  \begin{cases}
    U = \frac{m e^V}{\int_\Omega e^V\dd x}, & x\in\Omega,\\
    0 = \Delta V - \mu + U, & x\in\Omega,\\
    \partial_\nu V = 0, & x\in\partial\Omega,\\
    \int_\Omega V\dd x = 0, \quad \mu = \frac{m}{\pi},
  \end{cases}
\end{equation}
under (sub)critical mass settings, i.e., $m\in(0,8\pi]$.

\begin{proposition}
  [\!\!\!{\cite[Theorem~2.2]{Horstmann2011}}]
  \label{prop: uniqueness of symmetry steady states}
  Let \(\Omega=B_1\subset\mathbb{R}^2\) and \(m\in(0,8\pi]\). 
  Suppose 
  the pair \((U,V)\in (C^0(\overline{\Omega})\cap C^2(\Omega))^2\) 
  of radially symmetric functions solves the system \eqref{sys: stationary problem of uv} classically.
  Then \((U,V) = (m/\pi, 0)\).   
\end{proposition}

We begin with the natural introduction of the problem~\eqref{sys: stationary problem of uv}, 
according to long-time asymptotic behaviors of solutions to the system~\eqref{sys: JL}, 
based on Proposition~\ref{prop: global boundedness}. 

\begin{proposition}
\label{prop: long-time behaviors}
Let \(\Omega = B_1\subset\mathbb{R}^2\). Suppose that \(u_0\in C^0(\overline{\Omega})\) is radially symmetric and nonnegative such that
\begin{equation*}
    m = \int_\Omega u_0\dd x \in (0, 8\pi].
\end{equation*}
Then for any unbounded sequence \(\{t_k>0: k\in\mathbb{N}\}\), 
up to the extraction of a subsequence,
there exists 
 a pair \((U,V)\in(C^2(\overline\Omega))^2\) of radially symmetric functions
such that \((U,V)\) solves \eqref{sys: stationary problem of uv},
that  
\begin{equation*}
(u(x, t_k), v(x, t_k))\to (U,V)
\quad\text{in }C^2(\Omega)\times C^2(\Omega),
\end{equation*}
and such that
\begin{equation*}
  \mathcal{F}(u_0) - \mathcal{F}(U) 
  = \int_0^\infty\int_\Omega u|\nabla (\ln u - v)|^2,
\end{equation*}
where \((u,v)\) is the  solution of the system \eqref{sys: JL},
uniquely given by Proposition~\ref{prop: local well-posedness},
and $\mathcal{F}$ is defined by \eqref{sym: free energy}.
\end{proposition}

\begin{proof}
  The proof is a minor modification of \cite[Lemma~3.1]{Winkler2010a}.
  Thanks to Proposition~\ref{prop: global boundedness}, 
  it follows from standard parabolic and elliptic regularity theories that 
  there exists $\sigma\in(0,1)$ such that
  \begin{equation*} 
  (u,v)\text{ is bounded in } 
  (C^{2+\sigma, 1+\sigma/2}(\overline{\Omega}\times[1,\infty)))^2.
  \end{equation*} 
  So the free energy functional \eqref{sym: free energy} is bounded along the trajectory of $(u,v)$. 
  Integrating \eqref{eq: differential equation of free energy} over $(0,\infty)$, we get 
  \begin{equation*}
  \int_0^\infty\int_\Omega u|\nabla(\ln u - v)|^2 < \infty.
  \end{equation*}
  Then for the given sequence $\{t_k\}_{k\in\mathbb{N}}$,
  we have 
  \begin{equation*}
    \int_0^1\int_\Omega u(x,t+t_k)|\nabla(\ln u(x,t+t_k) - v(x,t+t_k))|^2
    \to 0 
    \quad \text{as } k\to\infty.
  \end{equation*}
  In addition, 
  possibly after extracting a subsequence, 
  we may assume by Arzel\`a-Ascoli theorem that 
  there exists  a pair $(U,V)$ of nonnegative and radially symmetric functions 
  such that 
  \begin{equation*}
    (u(\cdot,t+t_k),v(\cdot, t+t_k))
    \to (U(\cdot,t),V(\cdot,t))
    \quad\text{in } 
    (C^{2,1}(\overline{\Omega}\times[0,1]))^2,
    \quad\text{as } k\to\infty.
  \end{equation*}
  which implies
  \begin{equation}
  \label{eq: mass of long-time limit U}
    \int_\Omega V(\cdot, 0) = \int_\Omega v(\cdot, t_k) = 0,
    \quad
    \int_\Omega U(\cdot, 0) = \int_\Omega u(\cdot, t_k) = m,
  \end{equation}
  \(
    \partial_\nu V(x,0) 
    = \partial_\nu v(x,t_k) = 0
    \text{ on }\partial\Omega, 
  \)
  as well as
  \begin{equation}
    \Delta V - \mu + U = 0
    \quad\text{in }
    \overline{\Omega}\times[0,1].
  \end{equation}
  
  Let $\mathcal{S} := \{(x,t)\in\Omega\times(0,1): U(x,t) > 0\}$. 
  Then $\mathcal{S}$ is a nonempty set by \eqref{eq: mass of long-time limit U}. 
  Let $\mathcal{C}$ be a connected component of $\mathcal S$. Then $\mathcal{C}$ is relatively open in $\Omega\times(0,1)$. 
  By Fatou's lemma, 
    \begin{equation*}
    \int_{\mathcal{S}} U|\nabla (\ln U-V)|^2 
    \leq \liminf_{k\to\infty}\int_{\mathcal{S}} u(x,t+t_k)|\nabla (\ln u(x,t+t_k) - v(x,t+t_k))|^2 = 0.
  \end{equation*}
  We have $|\nabla (\ln U - V)| = 0$ over $\mathcal C$, 
  which implies $\ln U - V$ is a constant in $\mathcal C$ 
and thus there exists $C>0$ such that 
  \begin{equation*} 
  U = C\exp(V) \geq C\exp(-\|V\|_{L^\infty(\Omega)}) > 0
  \quad\text{on } \overline{\mathcal C}.
  \end{equation*}
  In particular, $\partial\mathcal C\cap\Omega\times(0,1) = \emptyset$, 
  and hence $\mathcal C$ is also a closed subset of $\Omega\times(0,1)$ with respect to the relative topology. 
  So we have $\mathcal{C} = \Omega\times(0,1)$ due to connectivity of $\Omega\times(0,1)$,
  and thus $U = C\exp(V) > 0$ over $\overline\Omega\times[0,1]$. 
  We can calculate by \eqref{eq: mass of long-time limit U}
  \begin{equation*}
    C = \frac{m}{\int_\Omega \exp(V)}.
  \end{equation*}
  and thereby complete the proof of Proposition~\ref{prop: long-time behaviors}.
\end{proof}

To deduce uniqueness of radial solutions to \eqref{sys: stationary problem of uv} under (sub)critical settings, previous works, 
e.g. \cite{Wang2019,Horstmann2011}, 
mainly focused on a nonlocal elliptic problem derived from \eqref{sys: stationary problem of uv} by substituting $U$ into the second equation. 
Instead, we make use of symmetry assumptions and insert $V$ into the first equation via \eqref{sym: mass distribution function} and thus treat the degenerate elliptic problem~\eqref{sys: stationary problem of scalar parabolic problem}.
To this end, we prepare some lemmas.

The following maximal-type auxiliary property was established in \cite[Lemma~3.4]{Winkler2019}.

\begin{lemma}
\label{le: maximal-type auxiliary lemma}
Let $\ow$ and $\uw$ be two functions from $C^1([0,1])\cap C^2((0, 1))$ such that 
$\uw_\xi \geq 0$ and $\ow_\xi \geq 0$ for all $\xi\in(0,1)$, 
that $\ow(1) = \uw(1)$, 
and such that
\begin{equation*}
\mathcal{Q}\uw \leq \mathcal{Q}\ow 
\quad\text{for all } \xi\in(0, 1),
\end{equation*}
where $\mathcal{Q}$ is given by \eqref{sym: operator Q}.
If moreover
\begin{equation*} 
\ow(\xi) \geq \uw(\xi)\quad \text{for all } \xi\in(0, 1), 
\text{ but } \ow \not \equiv \uw,
\end{equation*}
then there exists $C>0$ such that
\begin{equation*}
\ow(\xi) \geq \uw(\xi) + C \xi(1-\xi) 
\quad\text{for all } 
\xi \in(0, 1).
\end{equation*}
\end{lemma}

We need the following lemma and refer to \cite[Lemma~2.3]{Mao2024} for a complete proof.

\begin{lemma}
  \label{le: f>=cx(1-x)} 
Let $I$ be an interval.
Suppose that the mapping $I\ni \lambda\mapsto f(\cdot,\lambda)\in C^1([0,1])$ is continuous and $f(0,\lambda)=f(1,\lambda)=0$ for all $\lambda\in I$. 
Then for $\lambda_0\in I$ and $\epsilon>0$, 
there exists $\delta>0$ such that
\begin{equation*}
|f(\xi,\lambda)-f(\xi,\lambda_0)|
\leq \epsilon\xi(1-\xi)
\end{equation*}
for all $\xi\in(0,1)$ and $\lambda\in(\lambda_0-\delta,\lambda_0+\delta)\cap I$.
\end{lemma}

We observe that the stationary problem \eqref{sys: stationary problem of scalar parabolic problem} also admits a family of subsolutions.

\begin{lemma}
  \label{le: a family subsolution}
  Let $m\in(0,8\pi]$ and $b>0$. 
  Define
  \begin{equation}
  \label{sym: subsolution}
  \uw_b(\xi) := \frac{m b\xi}{b+1-\xi} 
  \quad \text{for } \xi\in(0,1).
  \end{equation}
  Then for each $b>0$, $\uw_b$ is convex and satisfies 
  \begin{equation*}
    \mathcal{Q}\uw_b < 0 \quad\text{for all } \xi \in (0,1).
  \end{equation*}
  Moreover,
  \begin{equation}
    \label{eq: asymptotics of uw}
    \lim_{b\searrow0}\uw_b = 0
    \quad\text{and}\quad
    \lim_{b\to\infty}\uw_b = m\xi
  \end{equation}
  for all $\xi\in(0,1)$.
\end{lemma}

\begin{proof}
Noting
\begin{align*}
  \uw_{b\xi} = \frac{m b}{b+1-\xi}
  + \frac{m b\xi}{(b+1-\xi)^2}
  = \frac{m(b+1)b}{(b+1-\xi)^2}
\end{align*}
and
\begin{equation*}
  \uw_{b\xi\xi} = \frac{2m(b+1)b}{(b+1-\xi)^3} > 0,
\end{equation*}
we compute 
\begin{align*}
  \mathcal{Q}\uw_b &=
  -\frac{8m\xi(b+1)b}{(b+1-\xi)^3}
  -\frac{m^2\xi(b+1)b^2}{\pi(b+1-\xi)^3}
  +\frac{m^2\xi(b+1)b}{\pi(b+1-\xi)^2}\\
  &=-\frac{m\xi(b+1)b}{\pi(b+1-\xi)^3}\cdot[8\pi + mb-m(b+1-\xi)]\\
  &=-\frac{m\xi^2(b+1)b}{\pi(b+1-\xi)^3} \cdot(8\pi-m+m\xi) < 0
\end{align*}
for all $\xi\in(0,1)$.
The identities in \eqref{eq: asymptotics of uw} are obvious.
\end{proof}

Now we are in a position to prove Proposition~\ref{prop: uniqueness of symmetry steady states}.

\begin{proof}
[Proof of Proposition~\ref{prop: uniqueness of symmetry steady states}]
Applying the mass distribution function~\eqref{sym: mass distribution function} 
to the solution $(U,V)$ of the stationary problem~\eqref{sys: stationary problem of uv}, 
we obtain that $W\in C^1([0,1])\cap C^3((0,1))$ solves the system~\eqref{sys: stationary problem of scalar parabolic problem} classically, 
and there exists $C>0$ satisfying 
\begin{equation*}
    C^{-1} 
    < W_\xi(\xi) 
    = \pi U(\sqrt{\xi}) 
    < C
    \quad\text{for all } \xi\in(0,1),
\end{equation*}
where
\begin{equation} 
  \label{sym: W}
  W(\xi) := 2\pi\int_0^{\sqrt{\xi}}U(r)r\dd r.
\end{equation}
Set 
\begin{equation*}
\mathcal{S} := \{ a > 0 : W \leq \ow_a \text{ for all } \xi\in(0,1)\},
\end{equation*}
where $\ow_a$ is given by \eqref{sym: supersolution of steady states problem}.
In light of the first asymptotic property in \eqref{eq: ow asymptotic},
Lemma~\ref{le: function with bound derivative can be bounded} ensures that there exists $\alpha>0$ such that $W \leq \ow_\alpha$ for all $\xi\in(0,1)$, 
which implies $\mathcal{S}$ is nonempty.
Since the mapping $(0,\infty)\ni a\mapsto \ow_a\in C^1([0,1])$ is clearly continuous, 
the set $\mathcal{S}$ is a closed subset of $(0,\infty)$ with respect to the relative topology.

We claim that the set $\mathcal{S}$ is also relatively open in $(0,\infty)$. 
Indeed, suppose that $\alpha\in\mathcal{S}$, 
then $W\leq \ow_\alpha$ for all $\xi\in(0,1)$.
Since $\ow_\alpha$ is a strictly supersolution of \eqref{sys: stationary problem of scalar parabolic problem} (see \eqref{eq: ow is a strict stationary supersolution}), 
$\ow_\alpha\not\equiv W$. Then we have by Lemma~\ref{le: maximal-type auxiliary lemma} that there exists $\varepsilon>0$ such that
\begin{equation*} 
\ow_\alpha \geq W + \varepsilon\xi(1-\xi)
\quad\text{for all } \xi\in(0,1).
\end{equation*}
Thanks to the continuity of the mapping $(0,\infty)\ni a\mapsto \ow_a\in C^1([0,1])$,
Lemma~\ref{le: f>=cx(1-x)} is applicable to existence of $\delta > 0$ with the property that
\begin{equation*}
\ow_a - \ow_\alpha  
\geq - \varepsilon\xi(1-\xi) 
\quad \text{for all } \xi\in(0,1),
\end{equation*}
for each $a\in(\alpha-\delta, \alpha + \delta)\cap(0,\infty)$.
Summing two estimates above up yields
\begin{equation*}
  \ow_a \geq W 
  \quad \text{for all } \xi\in(0,1) \text{ and }
  a\in(\alpha-\delta, \alpha + \delta)\cap(0,\infty),
\end{equation*}
which implies $(\alpha-\delta, \alpha + \delta)\cap(0,\infty)\in\mathcal{S}$.

Now we have $\mathcal{S} = (0,\infty)$ due to connectivity of $(0,\infty)$. 
Recalling the second asymptotic property in \eqref{eq: ow asymptotic},
we obtain that
\begin{equation}
  \label{eq: W <= mxi}
W \leq \lim_{a\to\infty}\ow_a 
= m\xi
\quad\text{for all } \xi\in(0,1).
\end{equation}
Define 
\begin{equation*}
\tilde{\mathcal{S}} := \{ b > 0: W\geq \uw_b \text{ for all } \xi\in(0,1)\},
\end{equation*}
with $\uw_b$ given by Lemma~\ref{le: a family subsolution}. 
After an analogous connectivity argument, we may assert that
$\tilde{\mathcal{S}} = (0,\infty)$ 
and thus
\begin{equation*} 
W\geq \lim_{b\to\infty}\uw_b 
= m\xi
\quad\text{for all } \xi\in(0,1), 
\end{equation*}
according to \eqref{eq: asymptotics of uw},
which implies $W\equiv m\xi$ by virtue of \eqref{eq: W <= mxi}. 
Equivalently, $U\equiv m/\pi$ in view of the definition of $W$ in \eqref{sym: W}.
\end{proof}

As a direct consequence of Proposition~\ref{prop: uniqueness of symmetry steady states} and Proposition~\ref{prop: long-time behaviors}, 
we immediately obtain the global asymptotic stability of spatial homogeneity with (sub)critical mass. 

\begin{corollary}
  \label{cor: global stability}
  Let \(\Omega=B_1\subset\mathbb{R}^2\). 
  Suppose that the initial function \(u_0\in C^0(\overline{\Omega})\) is nonnegative and radially symmetric such that 
  \[
  \int_\Omega u_0\dd x \in(0, 8\pi].
  \]
  Then the classical solution of the system~\eqref{sys: JL} is globally bounded and 
  \begin{equation*}
    u \to \fint_\Omega u_0\dd x \quad\text{in } C^2(\Omega),
  \end{equation*}
  as $t\to\infty$.
\end{corollary}

We actually establish a variant of logarithmic Hardy-Littlewood-Sobolev inequality \cite{Carlen1992}.

\begin{corollary}
Let \(\Omega=B_1\subset\mathbb{R}^2\).
Suppose that \(U\in C^0(\overline\Omega)\) is a nonnegative and radially symmetric function 
with (sub)critical mass
\begin{equation*}
  \lambda := \int_\Omega U\in(0,8\pi].
\end{equation*}
Then
\begin{equation*}
  \mathcal{F}(U) = \int_\Omega U\ln U 
  - \frac{1}{2}\int_\Omega \bigg(U-\frac{\lambda}{\pi}\bigg)(-\Delta)^{-1}\bigg(U-\frac{\lambda}{\pi}\bigg)
  \geq \mathcal{F}\bigg(\frac\lambda\pi\bigg) 
  = \lambda \ln\frac{\lambda}{\pi},
\end{equation*}
where the equality \(\mathcal{F}(U) = \lambda\ln\lambda - \lambda\ln\pi\) holds if and only if
$U\equiv\lambda/\pi$. 
%Here, the inverse of minus Laplacian operator $(-\Delta)^{-1}(U-\fint_\Omega U\dd x)$ denotes the unique solution in $\mathcal{H}$ of Poisson equation $-\Delta v = U-\fint_\Omega U\dd x$ with homogeneous Neumann boundary conditions.
\end{corollary}

\section{Exponential convergence. Proof of Theorem~\ref{thm: global wellposedness}}

This section is devoted to exponential convergence to semi-trivial steady states of solutions to \eqref{sys: JL}, which completes the proof of Theorem~\ref{thm: global wellposedness}.

We begin with exponential convergence of solutions to \eqref{sys: scalar system derived from JL system} emanating from stationary supersolutions~\eqref{sym: supersolution of steady states problem} in the space of integrable functions. 

\begin{lemma}
  \label{le: Ma convergence rate}
  Let $\Omega =B_1\subset\mathbb{R}^2$.
  Then there exists a positive constant $C$ such that 
  whenever $m\in(0,8\pi]$ and $a>0$,
    \begin{equation}
    \label{eq: Ma >= mxi}
    \Ma \geq m\xi \quad \text{for all } (\xi,t)\in[0,1]\times[0,\infty),
  \end{equation}
  and 
\begin{equation}
  \label{eq: Ma - mxi exponentially decay}
  \sup_{t>0}\bigg\{\e^{Ct}\int_0^1(\Ma - m\xi)\dd\xi\bigg\} < \infty.
\end{equation}
  where the mass distribution function~$\Ma$ denotes the global classical solution to \eqref{sys: scalar system derived from JL system} corresponding to the initial datum $\ow_a$ defined by~\eqref{sym: supersolution of steady states problem}, 
as given in Lemma~\ref{le: mass distribution function}.
\end{lemma}

\begin{proof}
  Corollary~\ref{cor: global stability} asserts that 
  \begin{equation}
    \label{eq: Ma to m xi}
    \Ma\to m\xi \quad \text{in } C^2([0,1]),
    \text{ as } t\to\infty.
  \end{equation}
  In light of 
  \begin{equation*}
    \Ma(\xi,0) = %2\pi\int_0^{\sqrt\xi}\ua_0r\dd r = 
    \ow_a \geq m\xi
    \quad \text{for all } \xi\in(0,1),
  \end{equation*}  
  Lemma~\ref{le: comparision principle} warrants \eqref{eq: Ma >= mxi}.
  Write $\Ra:= \Ma - m\xi$ for $(\xi,t)\in[0,1]\times[0,\infty)$. 
  Then the function $\Ra$ is nonnegative according to \eqref{eq: Ma >= mxi} and satisfies 
  \begin{align}
    \label{eq: Rat}
    \Ra_t &= 4\xi \Ma_{\xi\xi} + \frac{\Ra\Ma_{\xi}}{\pi} \notag\\
    &= \bigg(4\xi \Ma_\xi - 4\Ma + \frac{{\Ra}^2}{2\pi}\bigg)_\xi + \frac{m\Ra}{\pi}\notag\\
    &= \bigg(4\xi \Ra_\xi - 4\Ra + \frac{{\Ra}^2}{2\pi}\bigg)_\xi + \frac{m\Ra}{\pi}
  \end{align} 
  for all $(\xi,t)\in(0,1)\times(0,\infty)$.
  Define 
  \begin{equation}
    \pa := \int_0^1\Ra\ln\frac{\Lambda}{\xi}\dd\xi \quad \text{for } t\in[0,\infty),
  \end{equation}
  where $\Lambda > 1$ is fixed such that 
  \begin{equation*} 
    \e > 2\Lambda,
  \end{equation*}
  thanks to
  \begin{equation*}
    \e = \sum_{i=0}^\infty \frac{1}{i!} > 1 + 1 = 2.
  \end{equation*}
  Multiplying \eqref{eq: Rat} by $\ln\Lambda - \ln\xi$ and integrating over $(0,1)$,
  we obtain 
  \begin{align}
    \label{eq: d/dt psia}
    \frac{\dd}{\dd t}\pa 
      &= 4\Ra_\xi(1,t)\ln\Lambda 
      + \int_0^1 \bigg(4\xi\Ra_\xi - 4\Ra + \frac{{\Ra}^2}{2\pi}\bigg)\frac{1}{\xi}\dd\xi 
      + \frac{m}{\pi}\int_0^1\Ra\ln\frac\Lambda\xi\dd\xi\notag\\
      &\leq - 4\int_0^1\frac\Ra\xi\dd\xi 
      + \int_0^1\frac{{\Ra}^2}{2\pi \xi}\dd\xi 
      + \frac{m}{\pi}\int_0^1\Ra\ln\frac{\Lambda}{\xi}\dd\xi
  \end{align}
  for all $t\in(0,\infty)$,
  where we use $\Ra_\xi(1,t)\leq0$ for all $t>0$, 
  since $\Ra$ is nonnegative and $\Ra(1,t)=0$ for all $t>0$.
  Noting 
  \begin{equation}
    \label{eq: exp x geq ex}
    \frac{1}{\xi} 
    = \Lambda^{-1}\exp\bigg(\ln\frac{\Lambda}{\xi}\bigg) 
    \geq \frac{\e}{\Lambda}\ln\frac{\Lambda}{\xi}
    \quad \text{for all } \xi > 0,
  \end{equation}
  we get from~\eqref{eq: d/dt psia}
  \begin{equation}
    \label{eq: d/dt psia -4e}
    \frac{\dd}{\dd t}\pa 
    \leq - 4\e^{-1}(\e-2\Lambda)\int_0^1\frac{\Ra}{\xi}\dd\xi 
    + \int_0^1\frac{{\Ra}^2}{2\pi\xi}\dd\xi 
  \end{equation}
  for all $t\in(0,\infty)$.
  By \eqref{eq: Ma to m xi}, we may find $\Ta \geq 0$ such that 
  \begin{equation*}
    \frac{\Ra}{2\pi} < 2\e^{-1}(\e-2\Lambda)
    \quad \text{for all } (\xi, t)\in(0,1)\times(\Ta, \infty),
  \end{equation*}
  and thus \eqref{eq: d/dt psia -4e} is reduced to 
  \begin{equation*}
    \frac{\dd}{\dd t}\pa 
    \leq - 2\e^{-1}(\e-2\Lambda)\int_0^1\frac{\Ra}{\xi}\dd\xi 
    \leq - 2\Lambda^{-1}(\e - 2\Lambda)\pa 
    \quad \text{for all } t > \Ta,
  \end{equation*}
  which implies 
  \begin{equation*}
    \pa \leq 4\pi\e^{-1}(\e-2\Lambda)\exp\big(-2\Lambda^{-1}(\e - 2\Lambda)(t-\Ta)\big)
    \quad \text{for all } t > \Ta.
  \end{equation*}
  This readily shows \eqref{eq: Ma - mxi exponentially decay} 
  with the constant $C$ given by 
  \begin{equation}
    \label{sym: C}
    C := 2\Lambda^{-1}(\e - 2\Lambda) > 0,
  \end{equation}
  because $\pa \geq \ln\Lambda \int_0^1(\Ma-m\xi)\dd\xi$ for all $t>0$.
\end{proof}

Analogously, we show exponential convergence of solutions to \eqref{sys: scalar system derived from JL system} emanating from stationary subsolutions~\eqref{sym: subsolution} in the space of integrable functions. 

\begin{lemma}
  \label{le: Mb convergence rate}
  Let $\Omega =B_1\subset\mathbb{R}^2$, 
  $m\leq8\pi$ and $b>0$.
  Denote $\Mb$ the mass distribution function of the solution to \eqref{sys: JL} corresponding to the initial datum $\uw_b$ defined by~\eqref{sym: subsolution}.
  Then 
  \begin{equation}
    \label{eq: Mb <= mxi}
    \Mb \leq m\xi \quad \text{for all } (\xi,t)\in[0,1]\times[0,\infty),
  \end{equation}
  and 
  \begin{equation}
  \label{eq: mxi - Mb exponentially decay}
  \int_0^1(m\xi - \Mb)\dd\xi 
  \leq \e^{-Ct}\bigg(\int_0^1(m\xi - \uw_b)\Big(1-\frac{\ln\xi}{\ln\Lambda}\Big)\dd\xi\bigg)
  \quad\text{for all } t>0,
\end{equation}
where the constant $C$ and $\Lambda\in(1,\e/2)$ are fixed constants taken from Lemma~\ref{le: Ma convergence rate}.
\end{lemma}

\begin{proof}
  Noting 
  \begin{equation*}
    \Mb(\xi,0) = \uw_b\leq m\xi 
    \quad\text{for all }\xi\in(0,1),
  \end{equation*}
  we see from Lemma~\ref{le: comparision principle} that \eqref{eq: Mb <= mxi} holds.
  Define 
  \begin{equation*}
    \pb := \int_0^1\Rb\ln\frac{\Lambda}{\xi}\dd\xi \quad \text{for } t\in[0,\infty),
  \end{equation*}
  where $\Rb:= m\xi - \Mb$ for $(\xi,t)\in(0,1)\times(0,\infty)$ and 
  $\Lambda > 1$ is taken from Lemma~\ref{le: Ma convergence rate}. 
  Then the nonnegative function $\Rb\geq0$ complies with 
  \begin{align}
    \label{eq: Rbt}
    \Rb_t &= - 4\xi \Mb_{\xi\xi} + \frac{\Rb\Mb_{\xi}}{\pi} \notag\\
    &= 4\xi \Rb_{\xi\xi} - \frac{\Rb \Rb_\xi}\pi + \frac{m\Rb}{\pi}\notag\\
    &= \bigg(4\xi \Rb_\xi - 4\Rb - \frac{{\Rb}^2}{2\pi}\bigg)_\xi 
    + \frac{m\Rb}{\pi}
  \end{align} 
  for all $(\xi,t)\in(0,1)\times(0,\infty)$.
  Multiplying \eqref{eq: Rbt} by $\ln\Lambda - \ln\xi$ and integrating over $(0,1)$,
  we obtain 
  \begin{align*}
    \frac{\dd}{\dd t}\pb
      &= 4\Rb_\xi(1,t)\ln\Lambda 
      + \int_0^1 \bigg(4\xi\Rb_\xi - 4\Rb - \frac{{\Rb}^2}{2\pi}\bigg)\frac{1}{\xi}\dd\xi 
      + \frac{m}{\pi}\int_0^1\Rb\ln\frac\Lambda\xi\dd\xi\\
      &\leq - 4\int_0^1\frac\Rb\xi\dd\xi 
      + \frac{m}{\pi}\int_0^1\Rb\ln\frac{\Lambda}{\xi}\dd\xi
  \end{align*}
  for all $t\in(0,\infty)$, because $\Rb_\xi(1,t)\leq0$ for all $t>0$.
  Recalling~\eqref{eq: exp x geq ex},
  we have 
  \begin{equation*}
    \frac{\dd}{\dd t}\pb 
    \leq - 4\Lambda^{-1}(\e-2\Lambda)\pb
    \quad\text{for all } t\in(0,\infty),
  \end{equation*}
  which implies 
  \begin{equation*}
    \pb \leq \exp\big(-4\Lambda^{-1}(\e - 2\Lambda)t\big)\bigg(\int_0^1(m\xi - \uw_b)\ln\frac{\Lambda}{\xi}\dd\xi\bigg)
    \quad \text{for all } t > 0.
  \end{equation*}
  As $\pb \geq \ln\Lambda \int_0^1(m\xi-\Mb)\dd\xi$ for all $t>0$,
  we verify \eqref{eq: mxi - Mb exponentially decay} 
  with the constant $C$ given by~\eqref{sym: C}.
\end{proof}

Now we are in a position to show Theorem~\ref{thm: global wellposedness}.

\begin{proof}
  [Proof of \upshape{Theorem~\ref{thm: global wellposedness}}]
  Let $M$ defined by~\eqref{sym: mass distribution function}, be the mass distribution function of $u$.
  Then one can find a constant $K>0$ such that 
  \begin{equation*}
    \frac{1}{K} < M_\xi(\xi,1) < K \quad \text{for all } \xi\in(0,1),
  \end{equation*}
  because $u(x,1)\in C^0(\overline{\Omega})$ is positive for all $x\in\overline{\Omega}$ by the strong maximum principle.
  In light of Lemma~\ref{le: function with bound derivative can be bounded},
  there exist constants $a, b\in(0,\infty)$ such that 
  \begin{equation*}
    \uw_b \leq M(\xi,1) \leq \ow_a 
    \quad \text{for all }\xi\in(0,1),
  \end{equation*}
  where $\ow_a$ and $\uw_b$ are defined in Lemma~\ref{le: a family of supersolutions of steady states} and Lemma~\ref{le: a family subsolution}, respectively.
  Thanks to Lemma~\ref{le: comparision principle}, 
  \begin{equation}
    \label{eq: Mb leq M leq Ma}
    \Mb(\xi,t) \leq M(\xi,t+1) \leq \Ma(\xi,t) 
    \quad \text{for all } (\xi,t)\in[0,1]\times[0,\infty),
  \end{equation}
  where $\Ma$ and $\Mb$ are taken from Lemma~\ref{le: Ma convergence rate} and Lemma~\ref{le: Mb convergence rate}, respectively.
  It follows from \eqref{eq: Mb leq M leq Ma} that 
  \begin{equation}
    \label{eq: estimate L1 of M-mxi}
    \int_0^1|M(\xi,t+1)-m\xi|\dd x 
    \leq \max\bigg\{\int_0^1(\Ma - m\xi)\dd\xi, \int_0^1(m\xi-\Mb)\dd\xi\bigg\}
  \end{equation}
  for all $t>0$.
  By virtue of Lemma~\ref{le: Ma convergence rate} and Lemma~\ref{le: Mb convergence rate}, 
  there exists a constant $C>0$ independent of $a,b\in(0,\infty)$ and $m\in(0,8\pi]$ such that 
  \begin{equation}
    \label{eq: int01 mxi - Mb}
    \sup_{t>0}\bigg\{\e^{Ct}\int_0^1(\Ma - m\xi)\dd\xi\bigg\} 
    + 
    \sup_{t>0}\bigg\{\e^{Ct}\int_0^1(m\xi - \Mb)\dd\xi\bigg\} 
    < \infty. 
  \end{equation}
  Recalling Gagliardo-Nirenberg-Sobolev interpolation inequality~\cite{Nirenberg1959}
\begin{equation*}
  \|f'\|_{L^{\infty}((0,1))} 
  \leq C_{\mathrm{GNS}} \|f\|_{W^{2,\infty}((0,1))}^{2/3}\|f\|_{L^1((0,1))}^{1/3} 
  \quad \text{for all } f\in W^{2,\infty}((0,1)),
\end{equation*}
  we have 
  \begin{align}
    \label{eq: u -fint u0 L infty}
    \pi\bigg\|u-\fint_\Omega u_0\dd x\bigg\|_{L^\infty(\Omega)} 
    &= \|(M-m\xi)_\xi\|_{L^\infty((0,1))}\notag\\
    &\leq C_{\mathrm{GNS}} \|M-m\xi\|_{W^{2,\infty}((0,1))}^{2/3}\|M-m\xi\|_{L^1((0,1))}^{1/3}
  \end{align}
  for all $t>1$.
  Since $u(\cdot,t+1)$ is uniform-in-$t$ bounded in $C^2(\overline{\Omega})$ by Corollary~\ref{cor: global stability}, 
we have $M(\cdot,t+1)$ is uniform-in-$t$ bounded in $C^2([0,1])$. 
Collecting \eqref{eq: estimate L1 of M-mxi}--\eqref{eq: u -fint u0 L infty}, 
we thereby obtain \eqref{eq: exponential stability} with $\ell:=C/3$.
\end{proof}

\section*{Acknowledgements}

The authors are grateful to the referees for their careful reading and many valuable comments.

The first author is funded by 
``the Fundamental Research Funds for the Central Universities'' (No.~B250201215)
and Basic Research Program of Jiangsu~(No.~BK20251482).
The second author has been supported by the Youth Research Fund of Anhui University of Technology (No. QZ202421, No. QD202377 ).
The third author is supported in part by National Natural Science Foundation of China (Nos. 12271092 and 11671079) and the Jiangsu Provincial Scientific Research Center of Applied Mathematics under Grant No. BK20233002.

\providecommand{\bysame}{\leavevmode\hbox to3em{\hrulefill}\thinspace}
\providecommand{\MR}{\relax\ifhmode\unskip\space\fi MR }
% \MRhref is called by the amsart/book/proc definition of \MR.
\providecommand{\MRhref}[2]{%
  \href{http://www.ams.org/mathscinet-getitem?mr=#1}{#2}
}
\providecommand{\href}[2]{#2}

\end{document}